\documentclass[12pt]{article}
\usepackage{graphicx}
\usepackage{amsfonts}
\usepackage{latexsym}
\usepackage{amsmath}
\usepackage{fancyhdr}

\makeatletter
\@addtoreset{figure}{section}
\def\thefigure{\thesection.\@arabic\c@figure}
\@addtoreset{table}{section}
\def\thetable{\thesection.\@arabic\c@table}

\def\@sect#1#2#3#4#5#6[#7]#8{\ifnum #2>\c@secnumdepth
     \def\@svsec{}\else
     \refstepcounter{#1}\edef\@svsec{\csname the#1\endcsname.\hskip .75em
}\fi
     \@tempskipa #5\relax
      \ifdim \@tempskipa>\z@
        \begingroup #6\relax
          \@hangfrom{\hskip #3\relax\@svsec}{\interlinepenalty \@M #8\par}%
        \endgroup
       \csname #1mark\endcsname{#7}\addcontentsline
         {toc}{#1}{\ifnum #2>\c@secnumdepth \else
                      \protect\numberline{\csname the#1\endcsname}\fi
                    #7}\else
        \def\@svsechd{#6\hskip #3\@svsec #8\csname #1mark\endcsname
                      {#7}\addcontentsline
                           {toc}{#1}{\ifnum #2>\c@secnumdepth \else
                             \protect\numberline{\csname the#1\endcsname}\fi
                       #7}}\fi
     \@xsect{#5}}
\def\@begintheorem#1#2{\it \trivlist \item[\hskip \labelsep{\bf #1\ #2.}]}
\def\section{\@startsection {section}{1}{\z@}{-3.5ex plus -1ex minus
 -.2ex}{2.3ex plus .2ex}{\normalsize\bf}}

\begin{document}

\title{The Structure of the Inverse to the Sylvester \\
Resultant Matrix}
\date{} 
\maketitle

\begin{center}
\author{Boris D. Lubachevsky\\
${ }$\\
{\em bdl@bell-labs.com}\\
${ }$\\
Bell Laboratories \\
600 Mountain Avenue \\
Murray Hill, New Jersey}
\end{center}

\setlength{\baselineskip}{0.995\baselineskip}
\normalsize
\vspace{0.5\baselineskip}
\vspace{1.5\baselineskip}

\begin{abstract}
Given polynomials $a ( \lambda )$ of degree
$m$ and $b( \lambda )$ of degree
$n$, we represent the inverse
to the Sylvester resultant matrix of
$a$ and $b$,
if this inverse exists, as a canonical sum of
$m~+~n$ dyadic matrices each of which
is a rational
function of zeros of $a$ and $b$.
As a result, we obtain the polynomial solutions
$x(\lambda)$ of degree $n-1$ and $y(\lambda)$ of degree $m-1$
to the equation
$a(\lambda) x(\lambda)  + b(\lambda) y(\lambda) = c(\lambda)$,
where $c(\lambda)$ is a given polynomial of degree $m+n-1$, as follows:
$x( \lambda )$ is the Lagrange interpolation polynomial
for the function $c( \lambda )/a( \lambda )$
over the set of zeros of $b( \lambda )$ and
$y( \lambda )$ is the one for the function $c( \lambda )/b( \lambda )$
over the set of zeros of $a( \lambda )$ .

${ }$\\
${ }$\\
${ }$\\
{\bf Key words}: interpolation polynomial, Lagrange, Hermite, 
fundamental polynomial, zero placement. single input, single output, 
adaptive control 


\end{abstract}

\[
\begin{array}{c}
{ }\\
{ }\\
{ }\\
{ }\\
{ }
\end{array}
\]

A method to solve the equation
\begin{equation}
\label{ax+by=c}
  a x + b y = c ,
\end{equation}
where $a = a ( \lambda )$, $b = b ( \lambda )$, and $c = c ( \lambda )$
are given, $x = x( \lambda )$, and $y = y( \lambda )$ are unknown
univariable polynomials of degrees $m, n, m+n-1, n-1$, and $m-1$, respectively,
is to solve the system of $m+n$ linear algebraic equations
\begin{equation}
\label{zS=d}
z S~=~d ,
\end{equation}
where $z ~=~z(x,y)$ and $d ~=~d(c)$ are $m+n$-dimensional row-vectors
composed of the coefficients
of polynomials $x, y$ and $c$, respectively,
$S~=~S(a,b)$ is the $(m+n) \times (m+n)$
$Sylvester$ $resultant$ matrix
of polynomials $a$ and $b$.
Solving equation \eqref{ax+by=c} 
constitutes an important single-input/single-output
case in the zero placement procedure
in control theory
(see, e.g., \cite{GoodwinSin}).
In this paper we explicitly represent the matrix $S ^{-1}$
as a canonical sum of $m+n$ dyadic matrices each of which
is a rational function of zeros of $a$ and $b$;
thus we give an explicit solution ($x( \lambda )$,$y( \lambda )$)
for \eqref{ax+by=c}.
This solution may have a practical application
in certain situations of adaptive control.
Although the formulation of this solution appears very simple
(see Corollaries\ 1 and 3),
the author is able to mention no other
work with this solution.

In the case in which
neither $a( \lambda )$ nor $b( \lambda )$ have multiple zeros,
the representation of $S^{-1}$ is much simpler than
in the general case.
First, we formulate the results for this special case.
\\

{\bf Theorem 1}.
\\
{\em Let}
$
~~a( \lambda ) = a_0
\lambda^m + a_1
\lambda^{m-1} +...+ a_m =
a_0 ( \lambda - \alpha_1 )
( \lambda - \alpha_2 )...( \lambda - \alpha_m )
$
\[
and~~~~~b( \lambda ) = b_0
\lambda^n + b_1
\lambda^{n-1} +...+ b_n =
b_0 ( \lambda - \beta_1 )
( \lambda - \beta_2 )...( \lambda - \beta_n )
~~~~~~~~~
~~~~~~~~~~~~~~~~~~~~~~~~~~
\]
{\em be complex polynomials,}
$a_0 b_0 \neq 0$,
{\em and let}
$S = S( a , b )$ 
{\em be their Sylvester matrix},
\[
S( a , b ) = \left [
\begin{array}{l l l l l l l l}
a_0 & a_1 & a_2 & . & . & .       & .       & 0 \\
0   & a_0 & a_1 & . & . & .       & .       & 0 \\
.   & .   & .   & . & . & .       & .       & 0 \\
0   & .   & .   & . & . & a_{m-1} & a_m     & 0 \\
0   & .   & .   & . & . & a_{m-2} & a_{m-1} & a_m \\
b_0 & b_1 & b_2 & . & . & .       & .       & 0 \\
0   & b_0 & b_1 & . & . & .       & .       & 0 \\
.   & .   & .   & . & . & .       & .       & 0 \\
0   & .   & .   & . & . & b_{n-1} & b_n     & 0 \\
0   & .   & .   & . & . & b_{n-2} & b_{n-1} & b_n 
\end{array}
\right ]
\begin{array} {c}
\uparrow \\
m \\
\mbox{rows} \\
  \\
\downarrow \\
\uparrow \\
n \\
\mbox{rows} \\
  \\
\downarrow
\end{array}
\]
\[
~~~\leftarrow ~~~~~~~~~n + m~~~~ \mbox{ columns} ~~~~~~~~~\rightarrow
\]
{\em If all zeros}
$\alpha_i$ {\em are simple
(i.e. pairwise different), and all zeros} $\beta_i$ 
{\em are simple, then
the adjoint matrix} $\mbox{adj}~S$ {\em is}
\newpage
\[ ~~~~\mbox{adj}~S(a,b) =
~~~~~~~~~~~~~~~~~~~~~~~~~~~~~~~~~~~~~~~~~~~~~~
~~~~~~~~~~~~~~~~~~~~~~~~~~~~~~~~~~~~~~~~~~~~~~
\]
\begin{equation}
\label{adjS}
~~~~~~~~~~~~~~~= ~a_0^n 
\sum_{i=1}^m
\left [
\begin{array}{c}
{\alpha_i}^{m+n-1}\\
{\alpha_i}^{m+n-2}\\
. \\
. \\
1
\end{array}
\right ]
\frac 
{
\begin{array}{c}
 \\
  \\
\prod \\
1 \le r \le m \\ 
r \neq i
\end{array}
b(\alpha_r)} 
{
 \displaystyle{
 \left.
 \frac {a(\lambda)} {\lambda-\alpha_i} 
 \right|_{\lambda = \alpha_i}
 }
}
\left [
0, ..., 0,\mbox{row}_m (\frac {a(\lambda)}{\lambda-\alpha_i})
\right ] +
\end{equation}
\[
+ (-1)^{mn}\, b_0^n 
\sum_{j=1}^n
\left [
\begin{array}{c}
{\beta_i}^{m+n-1}\\
{\beta_i}^{m+n-2}\\
. \\
. \\
1
\end{array}
\right ]
\frac 
{
\begin{array}{c}
 \\
  \\
\prod \\
1 \le s \le n \\ 
s \neq j
\end{array}
a(\beta_s)} 
{
 \displaystyle{
 \left.
 \frac {b(\lambda)} {\lambda-\beta_j}
 \right|_{\lambda = \beta_j}
 }
}
\left [
\mbox{row}_n (\frac {b(\lambda)}{\lambda-\beta_j}),0, ..., 0
\right ],
\]
{\em where} $\mbox{row}_k(p(\lambda))$ {\em denotes row}
$[ p_0 , p_1 , ... p_{k-1} ]$
{\em composed of the coefficients of the polynomial}
$p(\lambda) = p_0 \lambda ^{k-1} + p_1 \lambda ^{k-2} + ... +p_{k-1}$.
\\

{\bf Corollary 1}.
\\
{\em In the assumptions of Theorem 1, if} $\mbox{det} S(a,b)  \neq 0$,
{\em then}
\[
S(a,b)^{-1} =
~~~~~~~~~~~~~~~~~~~~~~~~~~~~~~~~~~~~~~~~~~~~~~~~
~~~~~~~~~~~~~~~~~~~~~~~~~~~~~~~~~~~~~~~~~~~~~~~~
\]
\begin{equation}
\label{S-1}
~~~~= ~ 
\sum_{i=1}^m
\left [
\begin{array}{c}
{\alpha_i}^{m+n-1}\\
{\alpha_i}^{m+n-2}\\
. \\
. \\
1
\end{array}
\right ]
\frac 
{1} 
{
\displaystyle{
 \left.
 b(\alpha_i) \left(
 \frac {a(\lambda)} {\lambda-\alpha_i}
 \right|_{\lambda = \alpha_i}
 \right)
}
}
\left [
0, ..., 0,\mbox{row}_m (\frac {a(\lambda)}{\lambda-\alpha_i})
\right ] +
\end{equation}
\[
~~~+ ~
\sum_{j=1}^n
\left [
\begin{array}{c}
{\beta_i}^{m+n-1}\\
{\beta_i}^{m+n-2}\\
. \\
. \\
1
\end{array}
\right ]
\frac 
{1} 
{
\displaystyle{
 \left.
 a(\beta_j) 
 \left(
  \frac {b(\lambda)} {\lambda-\beta_j}
 \right|_{\lambda = \beta_j}
 \right)
}
}
\left [
\mbox{row}_n (\frac {b(\lambda)}{\lambda-\beta_j}),0, ..., 0
\right ],
\]
{\em and the solution to equation} \eqref{ax+by=c} 
{\em is given by the formulas}
\begin{equation}
\label{laginter}
x(\lambda)~=~ 
\sum_{j=1}^n
\frac {
 \displaystyle{
 \frac {b(\lambda)} {\lambda - \beta_j}
 }
}
{
 \displaystyle{
 \left.
 \frac {b(\lambda)} {\lambda - \beta_j} 
 \right|_{\lambda = \beta_j}
 }
}
~\frac {c(\beta_j)} {a(\beta_j)},
~~~~~y(\lambda)~=~
\sum_{i=1}^m
\frac {
 \displaystyle{
 \frac {a(\lambda)} {\lambda - \alpha_i}
 }
}
{
 \displaystyle{
 \left.
 \frac {a(\lambda)} {\lambda - \alpha_i} 
 \right|_{\lambda = \alpha_i}
 }
}
~\frac {c(\alpha_i)} {b(\alpha_i)},
\end{equation}
{\em i.e.,} 
$x( \lambda )$ 
{\em is the Lagrange interpolation polynomial}
{\em for the function }
$c( \lambda )/a( \lambda )$
{\em over the set of zeros of}
$b( \lambda )$ 
{\em and}
$y( \lambda )$ 
{\em is the one for the function}
$c( \lambda )/b( \lambda )$
{\em over the set of zeros of}
$a( \lambda )$ .
\\

The following corollary describes the
asymptotic structure of $S^{-1}$,
$x( \lambda )$, and
$y( \lambda )$ when a zero of
$a( \lambda )$ approaches a zero
of $b( \lambda )$.
\\

{\bf Corollary 2}.
\\
{\em Let the coefficients of the polynomials}
$a( \lambda ) = a_{\tau} ( \lambda )$ 
{\em and}
$b( \lambda ) = b_{\tau} ( \lambda )$
{\em depend on a parameter} 
$\tau$, $\tau \in \{ \tau \}$.
{\em Denote by}
$x_{\tau} ( \lambda )$ and $y_{\tau} ( \lambda )$ 
{\em the corresponding
solutions of} \eqref{ax+by=c}.
{\em Let}
$\tau_*$ {\em be an accumulation point in} \{$\tau$\},
{\em and assume that for all}
$\tau  \neq  \tau_*$ 
{\em the zeros} 
$\alpha_i = \alpha_i ( \tau )$
{\em and} 
$\beta_j = \beta_j ( \tau )$
{\em are simple and that for} 
$\tau  \rightarrow  \tau_*$ 
{\em we have}
\[
a_{\tau} ( \lambda )  \rightarrow  a_* ( \lambda ),~
\alpha_i ( \tau ) \rightarrow {\alpha_i}^*,~ i = 1,...,m,~~~
b_{\tau} ( \lambda )  \rightarrow  b_* ( \lambda ),~
\beta_j ( \tau ) \rightarrow {\beta_j}^*,~ j = 1,...,n.
\]
{\em If}
\[ {\alpha_1}^* = {\beta_1}^* 
\stackrel{\rm def}{=}
 \theta
\]
{\em and}
${\alpha_i}^* \neq {\beta_j}^*$ 
{\em for all pairs}
$(i, j)$, 
{\em different from}
($i=1, j=1$)
{\em (zeros of}
$a_* ( \lambda )$ 
{\em and}
$b_* ( \lambda )$
{\em are not necessarily simple), then}
\[
\begin{array}{c}
{ }\\
\mbox{lim}\\
\tau \rightarrow {\tau}_* 
\end{array}
(\alpha_1(\tau) - \beta_1(\tau)) 
S{(a_{\tau},b_{\tau})}^{-1} =
\left [
\begin{array}{c}
{\theta}^{m+n-1}\\
{\theta}^{m+n-2}\\
.\\
.\\
1
\end{array}
\right ]
\frac {1} {
\left.
\displaystyle{
 \frac {a_* (\lambda) b_* (\lambda) } {(\lambda-\theta)^2}
}
\right|_{\lambda=\theta}
}
\left [
-\mbox{row}_n
(\frac{b_*(\lambda)}{\lambda-\theta})
,~
\mbox{row}_m
(\frac{a_*(\lambda)}{\lambda-\theta})
\right ],
\]

\[
\begin{array}{c}
{}\\
\mbox{lim}\\
\tau \rightarrow {\tau}_*
\end{array}
(\alpha_1(\tau) - \beta_1(\tau))
x_{\tau} (\lambda) =
-\, \frac {c(\theta)} 
{
\left.
 \displaystyle{
 \frac{a_*(\lambda)b_*(\lambda)}{(\lambda-\theta)^2}
 }
\right|_{\lambda = \theta}}
~\frac {b_* (\lambda)}{\lambda-\theta} \,,
\]

\[
\begin{array}{c}
{}\\
\mbox{lim}\\
\tau \rightarrow {\tau}_*
\end{array}
(\alpha_1(\tau) - \beta_1(\tau))
y_{\tau} (\lambda) =
\frac {c(\theta)} 
{\left.
\displaystyle{
 \frac{a_*(\lambda)b_*(\lambda)}{(\lambda-\theta)^2}
}
\right|_{\lambda = \theta}}
~\frac {a_* (\lambda)}{\lambda-\theta} \,.
\]
\\

In the case of simple zeros the representation of $S^{-1}$
relates to the Lagrange interpolation formula.
In the general case
the representation of $S^{-1}$ shall relate to
the general Hermite polynomial interpolation formula.
A version of this formula is presented below for reference.
\\
\\
\\
Let
\begin{equation}
\label{a-decomp}
a( \lambda ) = a_0 ( \lambda - \alpha_1 )^{m_1}
( \lambda - \alpha_2 )^{m_2}
...
( \lambda - \alpha_s )^{m_s}
\end{equation}
be a complex polynomial of degree $m = m_1 + m_2 +...+m_s$,
with zeros $\alpha_1 ,..., \alpha_s$,
\ \ $m_k>0, k=1,...,s$,
\ \ $a_0 \neq 0$,
\ \ $\alpha_k \neq \alpha_i$  for $k \neq i$.
Let $f( \lambda )$ be a complex function
which for $k = 1, 2, ..., s$ is defined at $\lambda = \alpha_k$ and
has $m_k -1$ successive derivatives
$f^{[1]} ( \alpha_k ),
f^{[2]} ( \alpha_k )
,...,
f^{[m_k -1]} ( \alpha_k )$.
The polynomial $p( \lambda )$ of degree no larger than $m-1$
for interpolation $f( \lambda )$ is as follows.
\[
p(\lambda) = 
\begin{array}{c}
m_1 -1 \\
\sum \\
i = 0
\end{array}
f^{[i]} (\alpha_1) u_{1,i} (\lambda ) +
\begin{array}{c}
m_2 -1 \\
\sum \\
i = 0
\end{array}
f^{[i]} (\alpha_2) u_{2,i} (\lambda ) + ... +
\begin{array}{c}
m_s -1 \\
\sum \\
i = 0
\end{array}
f^{[i]} (\alpha_s) u_{s,i} (\lambda ) 
\]
where $u_{k,i} ( \lambda )$,
$i = 0,1,...,m_k -1$,
$k = 1,2,...,s$,
are $m$ polynomials of degrees no larger than $m-1$.
For argument values $\lambda = \alpha_k$, $k=1,...,s$,
the values of the polynomial $p( \lambda )$ 
and its $m_k -1$ successive
derivatives
are the same as those of the function $f( \lambda )$.

The polynomials $u_{k,i}$ are independent of $f$
and are uniquely defined given
the zeros of the polynomial $a(\lambda)$ in \eqref{a-decomp} 
with their multiplicities.
They are called the {\em fundamental} polynomials.
For any particular set
$\{ m_1 ,m_2 ,...m_s \}$
explicit expressions for the fundamental polynomials
can be derived \cite{LanTis,Ralston}.
In the statement of Theorem\ 2 below we assume availability of
the fundamental polynomials $u_{k,i}$ corresponding to the polynomial $a$
as in \eqref{a-decomp}.

Similarly, $n$ fundamental polynomials $v_{\ell,j} ( \lambda )$
of degree no larger than $n-1$
correspond to a complex polynomial $b( \lambda )$ of degree $n$ as follows.
If
\begin{equation}
\label{b-decomp}
b( \lambda ) = b_0 ( \lambda - \beta_1 )^{n_1}
( \lambda - \beta_2 )^{n_2}
...
( \lambda - \beta_t )^{n_t},
\end{equation}
the zeros of $\beta ( \lambda )$ are $\beta_1 ,..., \beta_t$,
$n_{\ell} >0,~\ell = 1,...,t$,
$\beta_{\ell} \neq \beta_j$ for $\ell \neq j$,
then
\[
p(\lambda) = 
\begin{array}{c}
n_1 -1 \\
\sum \\
j = 0
\end{array}
f^{[j]} (\beta_1) v_{1,j} (\lambda ) +
\begin{array}{c}
n_2 -1 \\
\sum \\
j = 0
\end{array}
f^{[j]} (\beta_2) v_{2,j} (\lambda ) + ... +
\begin{array}{c}
n_t -1 \\
\sum \\
j = 0
\end{array}
f^{[i]} (\beta_t) v_{t,j} (\lambda ) 
\]
is the corresponding interpolation polynomial.
\\

We introduce the notation $\Lambda_k ( \lambda )$
for the following $\lambda$-column of height $k$
\[
\Lambda_k 
\stackrel{\rm def}{=}
\left [
\begin{array}{c}
\lambda^{k-1}\\
\lambda^{k-2}\\
. \\
. \\
\lambda \\
1
\end{array}
\right ]
\]
If $k = m+n$, the subscript can be omitted,
$\Lambda ( \lambda )  
\stackrel{\rm def}{=}  
\Lambda_{m+n} ( \lambda )$.
Theorem\ 1, in particular, implies that
in the case of simple zeros
each column of $S^{-1}$ is
a linear combination of
the columns obtained as the values
of $\Lambda ( \lambda )$
when $\lambda$ takes on the values of the zeros
of the polynomials $a$ and $b$.
The following Theorem\ 2, in particular,
asserts that in the general case,
each column of $S^{-1}$ is a linear combination
of the columns obtained as the values
of $\Lambda ( \lambda )$ and its derivatives
\[
\Lambda^{[1]} (\lambda)
=
\left [
\begin{array}{c}
(m+n-1)\lambda^{m+n-2}\\
(m+n-2)\lambda^{m+n-3}\\
. \\
. \\
2\lambda \\
1 \\
0
\end{array}
\right ] ,~~
\Lambda^{[2]} (\lambda)
=
\left [
\begin{array}{c}
(m+n-1)(m+n-2)\lambda^{m+n-3}\\
(m+n-2)(m+n-3)\lambda^{m+n-4}\\
. \\
. \\
2 \\
0 \\
0
\end{array}
\right ],
...
\]
when $\lambda$ takes on the values of the zeros
of the polynomials $a$ and $b$.
\\

{\bf Theorem 2}.
\\
{\em Let}
$a( \lambda )$ 
{\em be a complex polynomial of degree}
$m$ {\em as in} \eqref{a-decomp}
{\em and let} $u_{k,i} ( \lambda )$, $k=1,...,s$, $i=1,...,m_k$,
{\em be the} $m$ {\em corresponding fundamental polynomials.
Let} $b(\lambda )$ {\em be a complex polynomial
of degree} $n$ {\em as in} \eqref{b-decomp},
{\em and let} $v_{\ell,j} (\lambda )$, $\ell=1,...,t$, $j=1,...,n_{\ell}$,
{\em be the} $n$ {\em corresponding fundamental polynomials.
Let} $S = S(a,b)$ {\em be the} $(m+n) \times (m+n)$ {\em Sylvester
matrix defined as in Theorem\ 1.}

{\bf I}. {\em If} $\mbox{det}S\neq0$,
{\em then there exist} $m$ {\em polynomials}
$U_{k,i}(\lambda)$, $k=1,...,s$, $i=1,...,m_k$,
{\em of degree no larger than} $m-1$,
{\em and} $n$ {\em polynomials}
$V_{\ell,j}(\lambda)$, $\ell=1,...,t$, $j=1,...,n_{\ell}$,
{\em of degree no larger than} $n-1$,
{\em such that}
\[
S(a,b)^{-1} =
~~~~~~~~~~~~~~~~~~~~~~~~~~~~~~~~~~~~~~~~~~~~~~~~~~~~~
~~~~~~~~~~~~~~~~~~~~~~~~~~~~~~~~~~~~~~~~~~~~~~~~~~~~~
\]
\begin{equation}
\label{S-1}
~~~~~~~~~~= \sum_{i=0}^{m_1 -1} 
\Lambda^{[i]} (\alpha_1) 
\left [
0, ... ,0, \mbox{row}_m(U_{1,i} (\lambda))
\right ] + ... +
\sum_{i=0}^{m_s -1}
\Lambda^{[i]} (\alpha_s)
\left [
0, ... ,0, \mbox{row}_m(U_{s,i} (\lambda))
\right ] 
~~~
\end{equation}
\[
~~~~~~~
+ \sum_{j=0}^{n_1 -1} 
\Lambda^{[j]} (\beta_1)
\left [
\mbox{row}_n(V_{1,j} (\lambda))
0, ... ,0
\right ] + ... +
\sum_{j=0}^{n_t -1}
\Lambda^{[j]} (\beta_t)
\left [
\mbox{row}_n(V_{t,j} (\lambda))
0, ... ,0
\right ] . 
~~~~
\]
{\em The polynomials} $U_{k,i}$
{\em and} $V_{\ell,j}$,
{\em which satisfy} \eqref{S-1}
{\em and the restrictions on their degrees,
are unique.}

{\bf II}.
{\em For each fixed} $k=1, 2,...,s$ 
{\em the} $m_k$ {\em polynomials} $U_{k,i}$,
$i=0,1,...,m_{k -1}$,
{\em can be found successively,
beginning with} $i=m_{k -1}$,
{\em by the following recurrence:}
\begin{equation}
\label{Urec}
U_{k,i} =
\frac{1}{b(\alpha_k)}
\left [
u_{k,i} (\lambda) - 
\sum_{i_1 = i+1}^{m_k -1}
\left (
\begin{array}{c}
i\\
i_1
\end{array}
\right )
b^{[i_1 - i]} (\alpha_k) U_{k,i_1} (\lambda)
\right ] ,
~ i = m_k -1, m_k -2, ..., 1, 0.
\end{equation}
{\em For each fixed} $\ell = 1,2,...,t$ 
{\em the} $n_\ell$ {\em polynomials} $V_{\ell,j}$,
$j=0,1,...,n_{\ell} -1$,
{\em can be found successively,
beginning with} $j = n_{\ell} -1$,
{\em by the following recurrence:}
\[
V_{\ell,j} =
\frac{1}{a(\beta_{\ell})}
\left [
v_{\ell,j} (\lambda) - 
\sum_{j_1 = j+1}^{n_{\ell} -1}
\left (
\begin{array}{c}
j\\
j_1
\end{array}
\right )
a^{[j_1 - j]} (\beta_{\ell}) V_{\ell,j_1} (\lambda)
\right ] ,
~ j = n_{\ell} -1, n_{\ell} -2, ..., 1, 0.
\]
\\

{\bf Example}.
\\
Let $\alpha_k$ be a zero of multiplicity $m_k =2$
of a polynomial $a( \lambda )$.
In the decomposition \eqref{S-1}
two dyadic summand matrices
correspond to this zero,
$\Lambda (\alpha_k ) [ 0,...,0, \mbox{row}_m (U_{k,0} (\lambda))]$ and
$\Lambda^{[1]} (\alpha_k) [ 0,...,0, \mbox{row}_m (U_{k,1}(\lambda)) ]$.
We will now derive an explicit expression for $U_{k,0}$ and $U_{k,1}$
using \eqref{Urec}.
Denote $\bar{a}_k ( \lambda ) = a(\lambda )/(\lambda - \alpha_k )^2$.
The corresponding fundamental polynomials are (see \cite{Ralston}):

\[
u_{k,0} ( \lambda )  = 
\frac {\bar{a}_k ( \lambda )} {\bar{a}_k ( \alpha_k )}
[
1 - \frac
{\bar{a}_k^{[1]} ( \alpha_k )( \lambda - \alpha_k )}
{\bar{a}_k ( \alpha_k )}
]
\]
\[
u_{k,1} ( \lambda )  = 
\frac
{\bar{a}_k ( \lambda )( \lambda - \alpha_k )} 
{\bar{a}_k ( \alpha_k )} .
~~~~~~~~~~~~~~~~~
\]
Following the procedure in Theorem\ 2,
we begin with the value $i  =  m_k -1  = 1$
in \eqref{Urec} and find
\[
U_{k,1} ( \lambda )  = 
\frac
{1} {b( \alpha_k )} u_{k,1} ( \lambda )  = 
\frac
{\bar{a}_k ( \lambda )( \lambda - \alpha_k )}
{\bar{a}_k ( \alpha_k )b( \alpha_k )}
\,.
\]
Note that the value of the sum 
$\displaystyle{\sum_{i_1 = i+1}^{m_k -1}}$
in \eqref{Urec} is assumed to be 0 in the above calculation,
since $m_k -1 < i+1$.
Next, for $i  =  m_k -2  = 0$
we calculate
\[
U_{k,0} ( \lambda )  =  \frac {1} {b( \alpha_k )}
[
u_{k,0} ( \lambda ) - b^{[1]} ( \alpha_k ) U_{k,1} ( \lambda )
]
= 
\]
\begin{equation}
\label{example}
~~~~~~~~~~~~~~ = 
\frac {\bar{a}_k ( \lambda )} {b( \alpha_k )\bar{a}_k ( \alpha_k )}
\left [
1 - 
( \lambda - \alpha_k )
\left (
\frac {\bar{a}_k^{[1]} ( \alpha_k )} {\bar{a}_k ( \alpha_k )}
-
\frac {b^{[1]} ( \alpha_k )} {b( \alpha_k )}
\right )
\right ]
\,.
\end{equation}
\\
   
{\bf Corollary\ 3}.
\\
{\em If} $\mbox{det} S(a,b) \neq 0$ 
{\em then Equation} \eqref{ax+by=c}
{\em has the following solution:}
$x( \lambda )$ 
{\em is the general Hermite interpolation polynomial 
for the function}
$c( \lambda )/a( \lambda )$ 
{\em over the set of zeros (with their multiplicities) of} 
$b( \lambda )$,
{\em and similarly}
$y( \lambda )$ 
{\em is the interpolation polynomial for the function}
$c( \lambda )/b( \lambda )$ 
{\em over the set of zeros (with their multiplicities) of}
$a( \lambda )$.
\newpage
{\bf Discussion}.
\\
If zeros of $a$ and $b$ are simple,
Theorem\ 1 represents $\mbox{adj} S$ in the form \eqref{adjS} .
This representation is valid for both cases,
$\mbox{det} S  =  0$ and $\mbox{det} S \neq 0$.
A representation for $\mbox{adj} S$
can be obtained if zeros are not simple
if we
multiply both sides of Equation \eqref{S-1}
by $\mbox{det} S$.
This representation is of a form similar to 
that in \eqref{S-1}: it is a linear combination 
of dyadic matrices.
The columns in the dyads are the same as in \eqref{S-1}.
The rows are produced by the operation ``row''
from
certain polynomials, 
which are equal to $U_{k,i}$ and $V_{\ell,j}$
up to a proportionality constant.
These ``row'' polynomials are unique
given the limitation on their degree
since $U_{k,i}$ and $V_{\ell,j}$
are unique polynomials.
However this general representation is
not necessarily valid in the case $\mbox{det} S  = 0$.
If zeros are simple, then
the zero terms
$a( \beta_j )$
and $b( \alpha_i )$
in the denominators
are canceled by
the corresponding terms in
$\mbox{det} S =  \prod_ {i=1}^m b( \alpha_i )
 =  \prod_ {j=1}^n a( \beta_j )$
which appear in the numerators.
If zeros are not simple
$\mbox{det} S$ does not necessarily
cancel all zero terms in the denominators.

In the example above 
it may be that
the coefficients of $b( \lambda )$ depend on a parameter $\tau$,
$ \tau \in \{ \tau \}$,
$b( \lambda ) = b_\tau ( \lambda )$,
the set $\{ \tau \}$ has an accumulation point $ \tau_*$,
and $b_\tau ( \alpha_k )  \rightarrow  0$
when $\tau \rightarrow \tau_*$,
but $b_{\tau} ( \alpha_k )~ \neq ~0$
for all $\tau ~ \neq ~ \tau_*$
and $| b_{\tau}^{[1]} ( \alpha_k )|  \ge  const >0$
for all $\tau \in \{ \tau \}$.
Because of the presence of the term
$\displaystyle{\frac {b^{[1]} ( \alpha_k )} {b ( \alpha_k )}}$
in expression \eqref{example},
$U_{k,0} ( \lambda ) \mbox{det} S$
retains $b( \alpha_k )$
in the denominator.
This,
in turn, entails that
when $ \tau   \rightarrow \tau_*$, the summand with
$\Lambda ( \alpha_k )$ tends to infinity
even after multiplication by $\mbox{det} S$.

Because of the uniqueness
of $U_{k,i}$ and $V_{\ell,j}$ in \eqref{S-1},
we can be assured that, 
in general,
no expression 
for $\mbox{adj} S$ of the form 
similar to \eqref{S-1}
exists which is valid
for both 
the cases $\mbox{det} S \neq 0$ and $\mbox{det} S = 0$.
\\

{\bf An application}.
\\
Solving the system of linear algebraic equations \eqref{zS=d} is easier
than finding the zeros of the polynomials $a( \lambda )$ and $b( \lambda )$.
Therefore formulas \eqref{laginter}(5) are inefficient 
for solving an isolated instance
of zero placement equation \eqref{ax+by=c}(1).

However, in the case in which the coefficients
of $a( \lambda )$ and $b( \lambda )$
are continuously changing with time $\tau$,
formulas \eqref{laginter} appear more attractive.
Instead of solving many instances of system \eqref{zS=d},
one may track zeros of $a( \lambda )$ and
$b( \lambda )$ as they change and substitute
these zeros into formulas \eqref{laginter}
or into their corresponding generalization.
This situation may arise if one
uses the zero placement procedure \eqref{ax+by=c}
in the course of adaptive control with adjustable
and time-variable $a( \lambda ) = a_\tau ( \lambda )$
and $b( \lambda ) = b_\tau ( \lambda )$.
\\

{\bf Proofs}.
\\
A proof
is only required for Theorem\ 2 and Corollary\ 3.
Whereas it is easy to establish Corollary\ 3, 
given the result of Theorem\ 2,
one can also establish it independently of Theorem\ 2,
if one computes 
$f(\lambda ) = a( \lambda )x( \lambda ) + b( \lambda )y( \lambda )$
and the appropriate number of derivatives of $f( \lambda )$
for each zero of $a$ and $b$.
\newpage
{\bf Proof of Theorem\ 2}.
\\
The notation $\mbox{row}_k (p( \lambda ))$ 
can be extended 
in an obvious way 
to the case of a column polynomial $p$
of degree no larger than $k-1$,
in which case $\mbox{row}_k (p( \lambda ))$
is a $k \times k$ matrix.
We have
\begin{equation}
\label{Arow=rowA}
A ( \mbox{row}_k (p( \lambda ))) =  \mbox{row}_k (Ap( \lambda ))
\end{equation}
\begin{equation}
\label{sumrow=rowsum}
\mbox{row}_k (p( \lambda )) + \mbox{row}_k (q( \lambda ))
=  \mbox{row}_k (p( \lambda )+q( \lambda ))
\end{equation}
for a $k \times k$ complex matrix $A$ and column $\lambda$-polynomials
$p( \lambda )$ and $q( \lambda )$ of degree no larger than $k-1$
and height $k$.

Let $Q$ be the $(m+n) \times (m+n)$ matrix in
the right-hand side of Equation \eqref{S-1}.
Applying rules \eqref{Arow=rowA} and \eqref{sumrow=rowsum} we obtain
\[
Q =  \mbox{row}_{m+n} (q( \lambda )),
\]
\[
q( \lambda ) = q_1 ( \lambda ) + q_2 ( \lambda ),
\]
where $q( \lambda ), q_1 ( \lambda ), q_2 ( \lambda )$
are column $\lambda$-polynomials
of degree no larger than $m+n-1$
and height $m+n$,
and the polynomials $q_1$ and $q_2$ are given by
\[
q_1 ( \lambda ) = 
\sum_{k=1}^s \sum_{i=0}^{m_k -1}
\Lambda^{[i]} ( \alpha_k )
U_{k,i} ( \lambda ),
~~~
q_2 ( \lambda ) = 
\lambda^m
\sum_{\ell=1}^t \sum_{j=0}^{n_{\ell} -1}
\Lambda^{[j]} ( \beta_{\ell} )
V_{\ell,j} ( \lambda ).
\]
The $(m+n) \times (m+n)$ identity matrix $I_{m+n}$
can be represented as
$I_{m+n}  =  \mbox{row}_{m+n} ( \Lambda_{m+n} ( \lambda ))$.
The condition $SQ = I_{m+n}$
can be rewritten as
\[
Sq( \lambda )  =  \Lambda ( \lambda ) .
\]
Using the specific structure of the Sylvester matrix $S$,
we have
\[
S \Lambda_{m+n} ( \lambda ) = 
\left [
\begin{array}{c}
\Lambda_n ( \lambda )a( \lambda )\\
\Lambda_m ( \lambda )b( \lambda )
\end{array}
\right ] .
\]
We can also calculate the corresponding $r$-th derivative
\[
S \Lambda^{[r]} ( \lambda )
= 
\left ( S \Lambda ( \lambda ) \right )^{[r]}  = 
\left [
\begin{array}{c}
\begin{array}{c}
r\\
\sum\\
i=0
\end{array}
\left ( 
\begin{array}{c}
i\\ 
r
\end{array}
\right )
\Lambda_n^{[i]} ( \lambda )a^{[r-i]} ( \lambda ) \\
\begin{array}{c}
r\\
\sum\\
i=0
\end{array}
\left ( 
\begin{array}{c}
i\\
r
\end{array}
\right )
\Lambda_m^{[i]} ( \lambda )b^{[r-i]} ( \lambda )
\end{array}
\right ] .
\]
Since 
$a( \alpha_k ) = a^{[1]} ( \alpha_k ) =
...=a^{[m_k -1]} ( \alpha_k ) = 0,~~k = 1,...,s$,
we have
\[
S \Lambda^{[r]} ( \alpha_k )  = 
\left [
\begin{array}{c}
0 \\
\begin{array}{c}
r\\
\sum\\
i=0
\end{array}
\left ( 
\begin{array}{c}
i\\ 
r 
\end{array}
\right )
\Lambda_m^{[i]} ( \alpha_k )
b^{[r-i]} ( \alpha_k )
\end{array}
\right ] ,~~
r = 0,1,...,m_k -1 .
\]
Summing all terms
$S \Lambda^{[r]} ( \alpha_k )U_{k,r} ( \lambda )$,
$r = 0,1,...,m_k -1$,
$k = 1,...,s$,
we obtain
\[
Sq_1 ( \lambda )  = 
\left [
\begin{array}{c}
0\\
h( \lambda )
\end{array}
\right ] ,
\]
where
\[
h( \lambda )  = 
\sum_{k=1}^ s
\sum_{r=0}^{m_k -1}
\sum_{i=0}^r
\left ( 
\begin{array}{c}
i\\ 
r
\end{array}
\right )
\Lambda_m^{[i]} ( \alpha_k )
b^{[r-i]} ( \alpha_k ) U_{k,r} ( \lambda )  = 
\]
\begin{equation}
\label{last-1}
~~~~~~~~~~~~~ = 
\sum_{k=1}^s
\sum_{i=0}^{m_k -1}
\Lambda_m^{[i]} ( \alpha_k )
\sum_{r=i}^{m_k -1}
\left ( 
\begin{array}{c}
i\\
r 
\end{array}
\right )
b^{[r-i]} ( \alpha_k ) U_{k,r} ( \lambda ) .
\end{equation}
According to \eqref{Urec}
the following identity holds:
\begin{equation}
\label{last}
u_{k,i} ( \lambda )  = 
\sum_{i_1 =i}^{m_k -1}
\left ( 
\begin{array}{c}
i\\
i_1
\end{array}
\right )
b^{[i_1 -i]} ( \alpha_k )U_{k,i_1} ( \lambda ) .
\end{equation}
The sum in the right-hand side of Equation \eqref{last}
is the same as the internal sum in \eqref{last-1}.
Thus, $h( \lambda )  =  \Lambda_m ( \lambda )$,
according to the general polynomial interpolation
formula and
\[Sq_1 ( \lambda )  =  
\left [ 
\begin{array}{c}
0\\
\Lambda_m ( \lambda )
\end{array}
\right ] \,.
\]

Similarly the identity
\[
Sq_2 ( \lambda )  =  
\left [ 
\begin{array}{c}
\lambda^m \Lambda_n ( \lambda ) \\
0 
\end{array}
\right ]
\]
can be established.
Hence representation \eqref{S-1} is proven.

The uniqueness of polynomials $U_{k,i}$ and $V_{\ell,j}$
easily follows from the non-singularity of the 
$(n+m) \times (n+m)$ {\em generalized Vandermonde matrix} 
\cite{LanTis}
\[
\left [
\Lambda ( \alpha_1 ),...,
\Lambda^{[m_1 -1]} ( \alpha_1 ),...,
\Lambda ( \alpha_s ),...,
\Lambda^{[m_s -1]} ( \alpha_s ),
\Lambda ( \beta_1 ),...,
\Lambda^{[n_1 -1]} ( \beta_1 ),...,
\Lambda ( \beta_t ),...,
\Lambda^{[n_t -1]} ( \beta_t )
\right ] .
\]
when $\alpha_k ~ \neq ~ \beta_{\ell}$,
for all pairs $k, \ell$,
$~k = 1,...,s$,$~\ell = 1,...,t$.
\hfill $\Box$
\\

{\bf  Acknowledgment}.
I would like to thank James Mc Kenna for a helpful comment
and Debasis Mitra for carefully reading the text.

\end{document}